\documentclass{amsproc}
\newif\ifpdf
    \ifx\pdfoutput\undefined
    \pdffalse % we are not running PDFLaTeX
    \else
    \pdfoutput=1 % we are running PDFLaTeX
    \pdftrue
    \fi

    \ifpdf
    \usepackage[pdftex]{graphicx}
    \else
    \usepackage{graphicx}
    \fi
\usepackage{amsmath,amssymb,amsthm,epsfig}

\newcommand\Ll{L_L}

\newcommand\rni{R_{NI}}

\newcommand\Z{\mathbb Z}

\newcommand\s{\mathcal S}

\newtheorem{theorem}{Theorem}[section]
\newtheorem{proposition}{Proposition}[section]
\newtheorem{lemma}[proposition]{Lemma}

\theoremstyle{definition}
\newtheorem{definition}[proposition]{Definition}

\theoremstyle{example}

\input{epsf.tex}

\begin{document}
\ifpdf
    \DeclareGraphicsExtensions{.pdf, .jpg, .tif}
    \else
    \DeclareGraphicsExtensions{.eps, .jpg}
    \fi

\title{Seesaw words in Thompson's group $F$}
\author{Sean Cleary}
\address{Department of Mathematics,
The City College of New York,
New York, NY 10031}
\email{cleary@sci.ccny.cuny.edu}
\thanks{The first author acknowledges support from PSC-CUNY grant \#64459-0033}
\author{Jennifer Taback}
\address{Department of Mathematics,
Bowdoin College,
Brunswick, ME 04011}
\email{jtaback@bowdoin.edu}
\thanks{The second author acknowledges support from
NSF grant DMS-0437481}


\begin{abstract}
We describe a family of elements in Thompson's group $F$ which
present a challenge to finding canonical minimal length
representatives for group elements, and which show  that $F$ is
not combable by geodesics. These elements have the property that
there are only two possible suffixes of long lengths for geodesic
paths to these elements from the identity; one is of the form $g^k$
and the other of the form $g^{-k}$ where $g$ is an element of a
finite generating set for the group.
\end{abstract}

\maketitle

\section{Introduction}
\label{sec:intro}

Thompson's group $F$ is a remarkable group with a poorly
understood but fascinating Cayley graph $\Gamma(F,\{x_0,x_1\})$,
using the standard finite presentation
$$F = \langle x_0,x_1 |
[x_0x_1^{-1},x_0^{-1}x_1x_0],[x_0x_1^{-1},x_0^{-2}x_1x_0^2]
\rangle.$$ We showed that this Cayley graph is not almost convex
\cite{ct}, meaning that metric balls are significantly folded in
upon themselves, that there are geodesics which cannot be extended
past a given point, ending in {\em dead end elements}, and that
these dead end elements are all of depth 2 \cite{ctcomb}.
 Belk and Bux show that this Cayley graph is
additionally not minimally almost convex \cite{belkbux}, and Guba
uses it to show that the Dehn function of the group is quadratic
\cite{gubaquad}. Burillo \cite{josegrowth}, Guba  \cite{gubagrowth}
and Belk and Brown \cite{belkbrown}
have considered the growth of the group by trying to estimate the
size of balls in this Cayley graph as well.

The goal of this paper is to exhibit a family of elements in $F$
which:

\begin{itemize}
\item show the difficulty of constructing canonical minimal length
representatives for elements of $F$, and \item exhibit the failure
of this Cayley graph to satisfy the $k$-fellow traveller property
for geodesics originating at the identity, and thus show that this
graph is not combable by geodesics.
\end{itemize}

We approach the problem of the existence of canonical minimal
length representatives for elements of $F$ in this generating set
by considering the conditions under which certain generators
decrease the word length of elements.  
 In a free group with a free generating set,
 there is always a unique generator which reduces the
word length of a nontrivial element $w$. In groups with relations, understanding
which generators reduce word length of given elements can give
insight into the geometry of the group.

Let $\Gamma$ be the Cayley graph of $F$ with respect to the finite
generating set $\{x_0,x_1\}$.  All word lengths in the arguments
below are computed with respect to this generating set.  Our
approach towards the construction of minimal length
representatives for group elements is to follow Fordham's outline
and begin with $w \in F$, then find a generator $g \in \{x_0^{\pm
1}, \ x_1^{\pm 1}\}$ so that $|wg| = |w| -1 $, where $| \cdot |$
denotes word length with respect to this generating set. Iteration
of this process will construct minimal length representatives, and
a natural goal is to find  a canonical way to proceed at elements
for which several generators decrease word length.

We describe below a class of `seesaw' words with the following
property.  If $w \in F$ is a seesaw word, then there is a unique
generator $g \in \{x_0^{\pm 1}, \ x_1^{\pm 1} \}$ so that $g$ and
$g^{-1}$ are the only two generators which decrease the word
length of $w$, that is, $|wg^{\pm 1}| < |w|$. Moreover, only $g$
decreases the word length of $wg^i$ for
 many iterations, and similarly only $g^{-1}$
decreases the word length of $wg^{-i}$ for many iterations.  The
existence of these words eliminates the possibility of a choice of
geodesics from the identity to each  elements satisfying the
$k$-fellow traveller property.

\section{Metric Properties of $F$}
\label{sec:F}

\begin{figure}\includegraphics[width=3.5in]{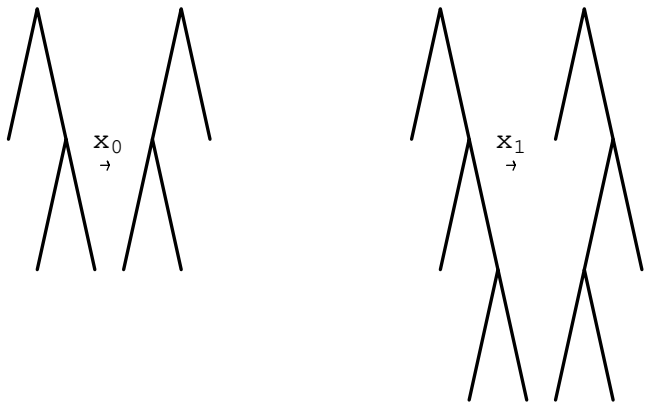}\\
\caption{The tree pair diagrams representing the generators $x_0$
and $x_1$ of $F$. \label{fig:x0x1}}
\end{figure}

We view elements of $F$ as pairs of finite rooted binary trees,
each with the same number of leaves.  To see the equivalence of
this with the group presentation, we refer the reader to Cannon,
Floyd and Parry \cite{cfp}.  We view our trees as consisting of a
collection of  `carets', which are interior nodes together with
their two downward-directed edges.  A leaf ending in a vertex of
valence one is called an {\em exposed leaf}.  A caret may have a
right child or a left child, if either or both of its leaves are
not exposed.  The tree pair diagrams representing the generators
$x_0$ and $x_1$ are given in Figure \ref{fig:x0x1}.

There is a natural reduction condition on tree pair diagrams to
ensure a unique tree pair diagram representing each group element.
Namely, a tree pair diagram $(T_-,T_+)$ is unreduced if there is a
caret with two exposed leaves, numbered $n$ and $n+1$ in both
trees of a tree pair diagram, and reduced otherwise.  To create a
reduced representative for an element, we simply remove these
common carets and renumber the exposed leaves.  We will assume
below that if $w=(T_-,T_+) \in F$, then the pair of trees is
reduced.  We refer to $T_-$ as the {\em negative tree} of the
pair, and $T_+$ as the {\em positive tree}.

There is an analytic interpretation of $F$ as a group of piecewise
linear homeomorphisms of the unit interval, subject to the
following two conditions:
\begin{enumerate} \item the slope of each linear piece is a power
of two, and \item the discontinuities of slope occur at points
whose coordinates are dyadic rationals.
\end{enumerate}
The carets in a binary rooted tree can be considered as
instructions for dyadic subdivision of the unit interval.  This
gives an equivalence between tree pair diagrams and the dyadic
piecewise linear homeomorphisms described above, where the trees
in the tree pair diagram are used to determine the domain and
range subdivisions for the homeomorphism.

With the analytic interpretation of $F$, group multiplication is
equivalent to composition of bijective functions.  In order to
multiply tree pair diagrams, one mimics the condition of bijective
function composition which requires that the range of one function
be the domain of the other.  Namely, to multiply $w = (T_-,T_+)$
and $v = (S_-,S_+)$, we create unreduced representatives of the
two elements, $(T'_-,T'_+)$ and $(S'_-,S'_+)$ respectively, in
which $T'_+ = S'_-$.  The product $wv$ is then represented by the
(not necessarily reduced) tree pair diagram $(T'_-,S'_+)$. See
Cannon, Floyd and Parry \cite{cfp} and Cleary and Taback
\cite{ctcomb} for details and examples of group multiplication of
tree pair diagrams.

\subsection{Fordham's method for computing word length in $F$}
Viewing elements of $F$ as tree pair diagrams, we use
Fordham's method \cite {blakegd} for computing the word length of
$w \in F$ with respect to the standard finite generating set
$\{x_0,x_1\}$ directly from the tree pair diagram representing
$w$. We now describe this remarkable method.

\begin{figure}
\includegraphics[width=3in]{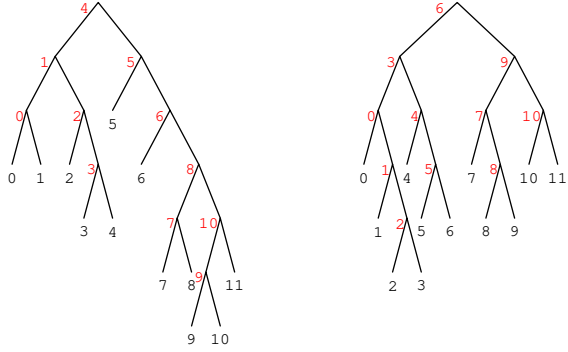}\\
\caption{The tree pair diagram for the group element
 $x_0^2 x_1 x_2 x_4 x_5 x_7 x_8 x_9^{-1} x_7^{-1} x_3^{-1} x_2^{-1} x_0^{-2}$
with carets and leaves numbered. \label{fig:treepairexample}}
\end{figure}

Fordham begins by dividing the carets in a binary rooted tree into
distinct types, roughly left, right and interior.  The left side
of the tree is the path of left edges beginning at the root caret;
the right side is defined analogously.  Left (respectively right) carets
have one leaf on the left (respectively right) side of the tree.  The
root caret of the tree is always considered a left caret. First,
the carets are numbered using an infix numbering scheme, beginning
with zero.  According to this infix order, the left child of a
caret is numbered before the caret, and any right child is
numbered after the caret.  Figure \ref{fig:treepairexample}
provides an example of a tree pair diagram with exposed leaves
numbered from left to right, and carets numbered in infix order in each tree.

Fordham's caret types are as follows:

\begin{enumerate}

\item $L_0$.  The first caret on the left side of the tree,
farthest away from the root caret. Every nonempty tree has exactly one
caret of type $L_0$, and it will always have infix number zero.

\item $\Ll$.  Any other left caret.

\item $I_0$. An interior caret which has no right child.

\item $I_R$. An interior caret which has a right child.

\item $R_I$. A right caret numbered $k$ with the property that
caret $k+1$ is an interior caret.

\item $R_0$. A right caret with no higher-numbered interior
carets.

\item $\rni$. A right caret which is neither an $R_I$ nor an $R_0$
caret.
\end{enumerate}
Let $w \in F$ be represented by the tree pair diagram $(T_-,T_+)$.
Fordham forms pairs of caret types by associating the carets in
$T_-$ and $T_+$ with the same caret number; for example, the first
pair consists of the type of caret zero in $T_-$ and the type of
caret zero in $T_+$ which will necessarily be $(L_0,L_0)$.
  Each pair of caret types is assigned a weight from the
following table.  Notice that the pattern of weights in the table
is symmetric around the diagonal.  The pair $(L_0,L_0)$ is
assigned weight $0$ and does not appear in the table.

\begin{center}
\begin{tabular}{|c|c|c|c|c|c|c|}
\hline
 & $R_0$ & $\rni$ & $R_I$ & $\Ll$ & $I_0$ & $I_R$ \\
 \hline
 $R_0$ & 0 & 2 & 2 & 1 & 1 & 3 \\ \hline
 $\rni$ & 2 & 2 & 2 & 1 & 1 & 3 \\ \hline
 $R_I$ & 2 & 2 & 2 & 1 & 3 & 3 \\ \hline
 $\Ll$ & 1 & 1 & 1 & 2 & 2 & 2 \\ \hline
 $I_0$ & 1 & 1 & 3 & 2 & 2 & 4 \\ \hline
$I_R$ & 3 & 3 & 3 & 2 & 4 & 4 \\ \hline
\end{tabular}
\end{center}

Fordham then proves the following theorem.

\begin{theorem}[Fordham \cite{blake}, Theorem 2.5.1]
\label{thm:blake} Given a element $w \in F$ described by the
reduced tree pair diagram $(T_-,T_+)$, the word length $|w|$ of
the element with respect to the generating set $\{x_0,x_1\}$ is
the sum of the weights of the caret pairings in $(T_-,T_+)$.
\end{theorem}

\subsection{Multiplication by the generators $x_0^{\pm 1}$ and
$x_1^{\pm 1}$}

Our method for finding minimal length representatives for elements
of $F$ with respect to the generators $\{x_0^{\pm 1},x_1^{\pm1}\}$
relies on being able to recognize changes in word length from
changes in caret types, using Fordham's method.  Fordham proves a
lemma which says that if the trees in a tree pair diagram $w =
(T_-,T_+)$ have the `correct' shape, then multiplication by a
generator will alter at most one pair of caret types in
$(T_-,T_+)$.  These conditions are seen to be necessary when one
tries to perform the multiplication by creating unreduced
representatives for the element $w$ and the generator.

\begin{lemma}[Fordham \cite{blakegd}, Lemma 2.3.1]
\label{lemma:conditions} Let $(T_-,T_+)$ be a reduced pair of
trees, each having $m+1$ carets, representing an element $x \in
F$, and $g$ a generator in  $\{x_0^{\pm 1}, \ x_1^{\pm
1}\}$.
\begin{enumerate}
\item If $g = x_0$, we require that the left subtree of the root
of $T_-$ is nonempty.

\item If $g = x_0^{-1}$, we require that the right subtree of the
root of $T_-$ is nonempty.

\item If $g = x_1$, we require that the left subtree of the right
child of the root of $T_-$ is nonempty.

\item If $g = x_1^{-1}$, we require that the right subtree of the
right child of the root of $T_-$ is nonempty.

\end{enumerate}
Then the reduced tree pair diagram for $x g$ also has $m+1$
carets, and there is exactly one $i$ with $0 < i \leq m$ so that
the pair of caret types of caret $i$ changes when $g$ is applied
to $x$.
\end{lemma}

Furthermore, Fordham shows that if the requirements of Lemma
\ref{lemma:conditions} are not met, then there is one additional
caret in the reduced tree pair diagram for $x g$ and $|xg| = |x| +
1$. The elements we describe in Section \ref{sec:seesaw} below
will be constructed to satisfy these conditions for $x_0$ and
$x_0^{-1}$.

It is easy to describe the exact change in a tree pair diagram $w
= (T_-,T_+)$ under right multiplication by one of the generators
$\{x_0^{\pm 1},x_1^{\pm1}\}$ when $w$ satisfies the conditions of
the above lemma. Each generator induces a rearrangement of the
subtrees of $T_-$, as can be seen either from performing
multiplication directly on the tree pair diagrams, or translating
back to the normal forms for elements.

\begin{figure}\includegraphics[width=4in]{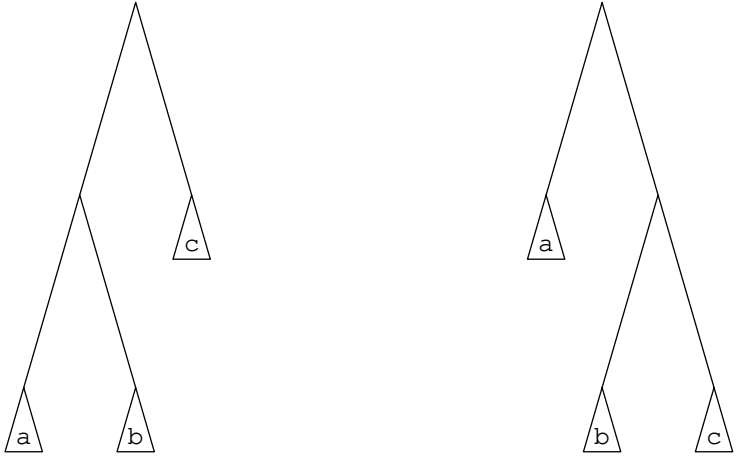}\\
\caption{This figure depicts from left to right the negative trees
from tree pair diagrams representing $w \in F$ and $wx_0$.
Multiplication by this generator performs a rearrangement of the
subtrees of this negative tree. \label{fig:x0rot}}
\end{figure}

\begin{figure}\includegraphics[width=4in]{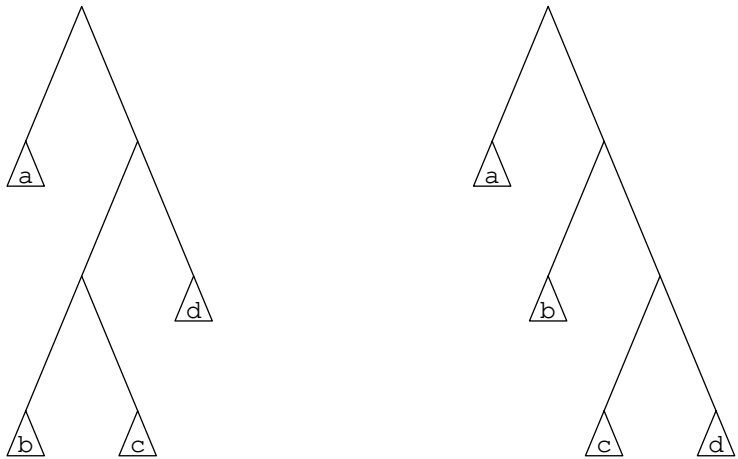}\\
\caption{ This figure depicts from left to right the negative
trees from tree pair diagrams representing $w \in F$ and $wx_1$.
Multiplication by this generator performs a rearrangement of the
subtrees of this negative tree.\label{fig:x1rot}}
\end{figure}

%Thus, the type of one of three carets will change under
%multiplication by these generators: the root caret, the right
%child of the root, or the left child of the root. We refer the
%reader to \cite{ctcomb} for an interpretation of these right
%multiplications as rotations of the negative tree at those
%ocations.

%We now
%describe the exact change in the tree pair diagram under right
%multiplication by a generator.
%
%
\begin{lemma}[\cite{ct}, Lemmas 2.6 and 2.7]\label{lemma:rotation}
\label{lemma:rotations} Let $g \in \{x_0^{\pm 1},x_1^{\pm1}\}$,
and $w = (T_-,T_+)$ an element of $F$ satisfying the condition of
Lemma \ref{lemma:conditions} corresponding to $g$. If $wg =
(S_-,S_+)$, then $S_+ = T_+$, and one of the following situations
applies to $S_-$.
\begin{enumerate}
\item If $g = x_0$ then $T_-$ is given by the left tree in Figure
\ref{fig:x0rot} and $S_-$ is given by the right tree in the
figure.
 \item If $g = x_0^{-1}$ then $T_-$ is given by the right tree in Figure
\ref{fig:x0rot} and $S_-$ is given by the left tree in the figure.
\item If $g = x_1$ then $T_-$ is given by the left tree in Figure
\ref{fig:x1rot} and $S_-$ is given by the right tree in the
figure. \item If $g = x_1^{-1}$ then $T_-$ is given by the right
tree in Figure \ref{fig:x1rot} and $S_-$ is given by the left tree
in the figure.
\end{enumerate}
\end{lemma}

Combining Lemmas \ref{lemma:conditions} and \ref{lemma:rotations},
we see that the single caret of $T_-$ which changes type under
multiplication by a generator is in one of three possible
positions: the root position, the right child of the root, or the
left child of the right child of the root.
%
%\begin{lemma}
%\label{lemma:newcaret} Let $w = (T_-,T_+)$ be an element of $F$,
%and let $g \in \{x_0^{\pm 1}, \ x_1^{\pm 1}\}$.  If $w$ does not
%satisfy the condition of lemma \ref{lemma:conditions}
%corresponding to $g$, then $|wg| \geq |w|$.
%\end{lemma}

\section{Background}
\label{sec:background}

\subsection{Minimal length representatives}

We now address the question of finding canonical minimal length
representatives for elements $w=(T_-,T_+) \in F$, when word length
is computed using the standard finite generating set.  In
\cite{ctcomb} we present a canonical method of constructing a
minimal length representative for a purely positive or purely
negative word, that is, one in which either $T_-$ or $T_+$ is
composed entirely of carets of type $R_0$.  We call this the
nested traversal method, as it creates a  minimal path based on
the order of the types of carets in the tree pair diagram.

The group $F$ also has an infinite presentation, namely
$$F=\langle x_n, \ n \geq 0 | x_i^{-1}x_jx_i = x_{j+1} \ \text{ if
}i<j \rangle.$$ There is a convenient set of normal forms for
elements of $F$ in this infinite presentation given by
$x_{i_1}^{r_1} x_{i_2}^{r_2}\ldots x_{i_k}^{r_k} x_{j_l}^{-s_l}
\ldots x_{j_2}^{-s_2} x_{j_1}^{-s_1} $ with $r_i, s_i >0$,
$i_1<i_2 \ldots < i_k$ and $j_1<j_2 \ldots < j_l$. This normal
form is unique if we further require that when both $x_i$ and
$x_i^{-1}$ occur, so does $x_{i+1}$ or  $x_{i+1}^{-1}$, as
discussed by Brown and Geoghegan in \cite{bg:thomp}. We note here
that replacing each occurrence of $x_n$ in the normal form of $w
\in F$ by $x_0^{-n+1}x_1x_0^{n-1}$ creates an expression for $w$
in terms of $\left\{x_0,x_1 \right\}$ which is usually not
minimal.

We attempt to create canonical minimal length representatives for
elements of $F$ by constructing a path of generators from an
element $w$ to the identity, each of which decreases word length
at the given point on the path. Crucial to this are the geometric
conditions detailed in Lemma \ref{lemma:conditions}.

Below we construct a family of words of $F$  we call {\em seesaw}
words, which have the following property.  If $w \in F$ is a
seesaw word, then there is a unique generator pair $g^{\pm 1}  \in
\{x_0^{\pm 1}, \ x_1^{\pm 1} \}$ which decreases the word length
of $w$.  Moreover, only multiplication by the generator $g$ (resp.
$g^{-1}$) decreases the word length of $wg$ (resp. $wg^{-1}$) for
many iterations.

\subsection{The $k$-fellow traveller property}

We now consider a collection of paths from the identity to each
vertex in a Cayley graph.  Such a collection of paths satisfies
the $k$-fellow traveller property if paths which end at points
distance one apart always stay within distance $k$ of each other.
This property is an important part of the definition of an {\em
automatic} group, and a {\em combable} group has canonical paths
from the identity to each element which satisfy the $k$-fellow
traveller property for some fixed $k$.  The elements we construct
below in Section \ref{sec:seesaw} show that no collection of
geodesic paths in the Cayley graph $\Gamma=\Gamma(F,\{x_0,x_1\})$
can satisfy this property.

All collections of paths we consider below in a given Cayley graph
include exactly one path for each group element.

Let $\gamma$ be a path in a Cayley graph $\Gamma(G,X)$ from the
identity to some element $w \in G$.  Then $\gamma = \Pi_{i=1}^n
g_i$ where $g_i \in X$.  We can view $\gamma$ as an eventually
constant map from $\Z^+$ into $G$ where $\gamma(i) =
\Pi_{j=1}^ig_j$ for $i \leq n$ and $\gamma(i) = w$ for $i>n$.  We
begin by defining the synchronous distance between two such paths.

\begin{definition}
Let $\gamma$ and $\eta$ be paths from the identity in the Cayley
graph $\Gamma(G,X)$ to elements $w$ and $v$, respectively.  Then
the {\em synchronous distance} between $\gamma$ and $\eta$ is
defined to be $$D_s(\gamma,\eta) = max_{i \in
\Z^+}d_{\Gamma}(\gamma(i),\eta(i)).$$
\end{definition}

We can now define the $k$-fellow traveller property for a pair of
paths; a collection of paths satisfies this property if every pair
of paths ending at vertices one unit apart satisfies the
$k$-fellow traveller property for the same constant $k$.

\begin{definition}
Two paths $\gamma$ and $\eta$ in the Cayley graph $\Gamma(G,X)$
from the identity to $w$ and $v$ respectively, with $d_{\Gamma}
(w,v) = 1$ are said to {\em $k$-fellow travel} if
$D_s(\gamma,\eta) \leq k$.
\end{definition}

\section{Seesaw words}
\label{sec:seesaw}

We now define the seesaw words mentioned above, and show that
Thompson's group $F$ contains arbitrarily large examples of such
words.

\begin{definition}
A element $w$ in a finitely generated group $G$ with finite
generating set $X$ is a {\em seesaw word of swing $k$ with respect
to a generator $g$} if the following conditions hold.  Let $|w|$
represent the word length of $w$ with respect to the generating
set $X$.
\begin{enumerate}
\item  Right multiplication by both $g$ and $g^{-1}$ reduces the
word length of $w$; that is, $|wg^{\pm 1}| = |w| - 1$, and for all
$h \in X  \smallsetminus \left\{g^{\pm1}\right\}$, we have
$|wh^{\pm 1}| \geq |w|$.

\item Additionally, $|wg^l| = |wg^{l-1}| - 1$   for integral $l
\in [1,k]$, and $|wg^{m}h^{\pm 1}| \geq |wg^{m}|$ for all $h
\in X  \smallsetminus \left\{g\right\}$ and integral $m \in
[1,k-1]$.

\item Similarly,  $|wg^{-l}| = |wg^{-l+1}| - 1$  for integral $l
\in [1,k]$, and $|wg^{-m}h^{\pm 1}| \geq |wg^{-m}|$ for all $h
\in X  \smallsetminus \left\{g^{-1}\right\}$  for integral $m \in
[1, k-1]$.
\end{enumerate}
\end{definition}

These are called seesaw words because they behave like a balanced
seesaw.  When in balance, there is a two-way choice about which
way to go down, but once that initial choice is made, there is
only the inexorable descent downward by the same generator for a
large number of steps determined by the swing.

Finite cyclic groups  ${\Z}_{2k}$ have seesaw words of swing $k$
with respect to the standard one-generator generating set at the
point $k$.  The only other examples of seesaw words of sizable
swing known to the authors beside those described here occur in
wreath products, such as $\Z \wr \Z$ and the lamplighter groups
$\Z_n \wr \Z$, as described in \cite{ctdeadlamp}.  All of those
wreath product examples are not finitely presentable.

\begin{theorem}
\label{thm:seesaw} Thompson's group $F$ contains seesaw words of
arbitrarily large swing  with respect to the generator $x_0$ in
the standard generating set $\{x_0, x_1\}$.
\end{theorem}

\begin{figure}\includegraphics[width=5.5in]{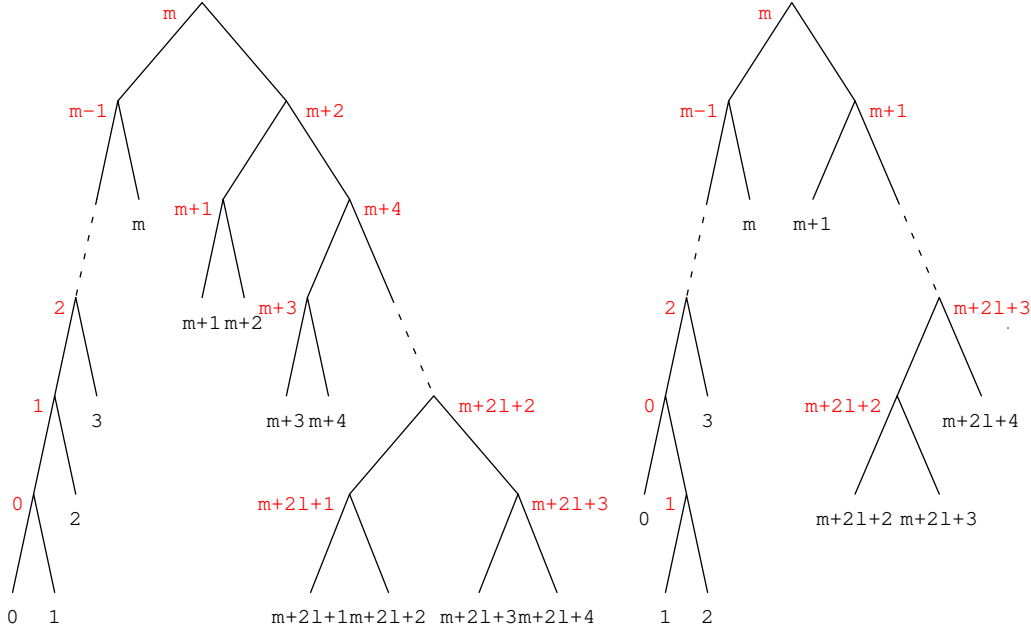}\\
\caption{The tree pair diagram representing the seesaw word
$x_0^{m-1} x_1 x_{m+2l+2} x_{m+2l+1}^{-1}
 x_{m+2l-1 }^{-1} \ldots  x_{m+3}^{-1} x_{m+1}^{-1} x_0^{-m}$.  Both trees have the same number of leaves.\label{seesawpic}}
\end{figure}

\begin{proof}
The idea of the proof is to construct elements $w$ with specific
pairs of caret types, chosen so that multiplication by both $x_0$
and $x_0^{-1}$ initially reduce the word length of $w$.  This is
easily seen using Fordham's methods.  Additionally, the pairs of
caret types which change under repeated multiplication by these
generators are also chosen so that these generators decrease word
length with each successive application.  One such family of words
is defined using two parameters, $l$ and $m$, and these words have
normal forms
$$x_0^{m-1} x_1 x_{m+2l+2} x_{m+2l+1}^{-1}
 x_{m+2l-1 }^{-1}  x_{m+2l-3}^{-1} \ldots  x_{m+3}^{-1} x_{m+1}^{-1} x_0^{-m}.$$
An example of a seesaw word of this form is given in Figure
\ref{seesawpic}.  We denote the family of these words by $\s$. The
parameter $l$ in the generic word of $\s$ given above determines
the length of the string of $R_I$ carets along the right side of
the negative tree of the pair, and of $R_{NI}$ carets on the right
side of the positive tree. The parameter $m$ determines the length
of the left sides of the trees. To ensure that our words in $\s$
have swing at least a given $k$, we let $l \geq k$ and  $m \geq
k$.

We consider what caret types the carets near the root of $T_-$ in
an element of ${\mathcal S}$ are paired with to ensure that
multiplication by both $x_0$ and $x_0^{-1}$ decrease word length
initially.  We see that for $w=(T_-,T_+) \in \s$, the root caret
of $T_-$, numbered $m$, is of type $L_L$ and is paired with caret
$m$ in $T_+$, also of type $L_L$.  Caret $m+2$, the right child of
the root in $T_-$, is of type $R_{I}$ and is paired with caret
$m+2$ in $T_+$, of type $R_{NI}$.

We now consider which generators reduce the length of $w$.  We
first note that all our words will satisfy the conditions of Lemma
\ref{lemma:conditions}, and thus the conclusions of Lemma
\ref{lemma:rotation} apply.

Right multiplication by $x_0$ will change caret $m$ in $T_-$ from
$L_L$ to $R_I$, so the pair of caret types will change from
$(L_L,L_L)$ to $(R_I,L_L)$, a reduction in weight from 2 to 1
which will reduce the overall word length by 1. Right
multiplication by $x_0^{-1}$ will change caret $m+2$ in $T_-$ from
$R_I$ to $L_L$, so the pair of caret types will change from
$(R_I,R_{NI})$ to $(L_L,R_{NI})$, a reduction in weight from 2 to
1 which will reduce the overall word length by 1.

Right multiplication by $x_1$ will change caret $m+1$ in $T_-$
from $I_0$ to $R_{NI}$, so the pair of caret types will change
from $(I_0, R_{NI} )$ to $(R_{NI}, R_{NI})$, an increase in weight
from 1 to 2 which will increase the overall word length by 1.

Right multiplication by $x_1^{-1}$ will change caret $m+2$  in
$T_-$ from $R_I$ to $I_R$, so the pair of caret types will change
from $(R_I, R_{NI} )$ to $(I_R, R_{NI})$, an increase in weight
from 2 to 3 which will increase the overall word length by 1.

Thus, $x_0$ and $x_0^{-1}$ reduce the word length of $w$ while
$x_1$ and $x_1^{-1}$ increase the word length.  Now we consider
how right-multiplication by each generator will affect the word
length of $w x_0^s$ for $s$ between $-l$ and $m$.

For $1 \leq  s  <  m$, the root caret of the negative tree of the
pair representing $wx_0^s$ will have caret number $m-s$ and be of
type $L_L$, and the right child of the root will be caret number
$m-s+1$ of type $R_{NI}$.   Both carets $m-s$ and $m-s+1$ will be
paired with carets of type $L_L$. These trees satisfy the
appropriate condition of Lemma \ref{lemma:conditions}.  When
$wx_0^s$ is multiplied by $x_0$, the pair of caret types which
changes corresponds to infix number $m-s$, and the change is from
$(L_L, L_L)$ to $(R_{NI},L_L)$ which will reduce length by 1.
Applying $x_0^{-1}$ will make the reverse change and increase
length.

Multiplication by $x_1^{-1}$ will change the pair of caret types
of caret $m-s+1$ either from $(R_{NI}, R_{NI})$ to $(I_0, R_{NI})$
or from $(R_{I}, R_{NI})$ to $(I_R, R_{NI})$, both of which
increase word length by one.  The trees representing $wx_0^s$ do
not satisfy the condition of Lemma \ref{lemma:conditions}
corresponding to the generator $x_1$, and thus it follows that
$|wx_0^sx_1| \geq |wx_0^s|$.

For $-l  <  s \leq -1$, the root caret of the negative tree of
$wx_0^s$ will be caret number $m-2s$ and have caret type $L_L$,
and the right child of the root will be caret number $m-2s+2$ of
type $R_{I}$. These will both be paired with carets of type
$R_{NI}$, so multiplication by $x_0^{-1}$ will change the pair of
caret types corresponding to carets $m-2s$ in both trees from
$(L_L, R_{NI})$ to $(R_{I},R_{NI})$, which will decrease word
length by 1. Multiplication by $x_0$ will make the reverse change
and increase word length by 1. Multiplication by $x_1^{-1}$ will
change the types of the carets of infix number $m-2s+2$ from from
$(R_I, R_{NI})$ to $(I_R, R_{NI})$ and also increase word length.
Multiplication by $x_1$ will change the types of the carets of
infix number $m-2s+1$ from $(I_0,R_{NI})$ to $(R_{NI},R_{NI})$,
increasing the word length by one.

Thus we see that all $w \in \s$ are seesaw words, and that there
are such words of any swing $k$.
\end{proof}

\begin{figure}\includegraphics[width=5in]{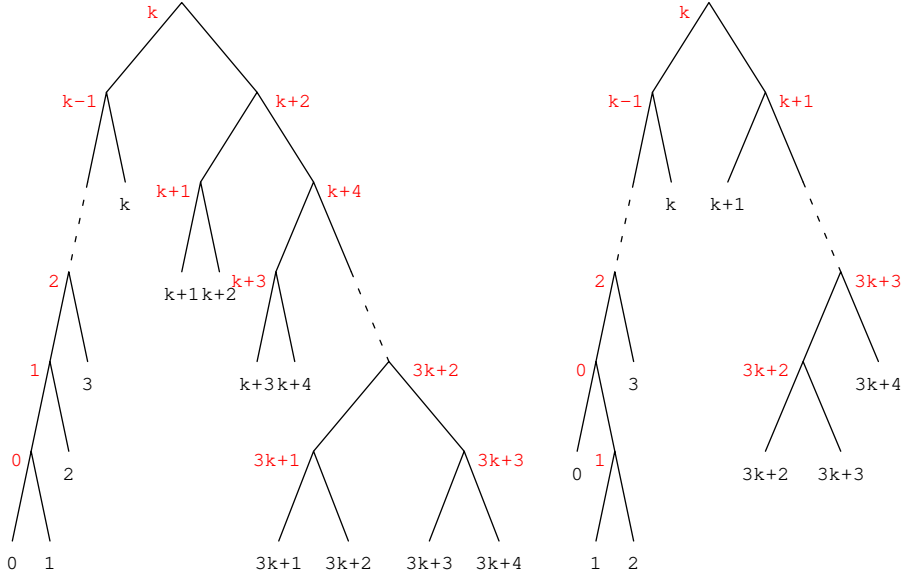}\\
\caption{Seesaw word  $x_0^{k-1} x_1 x_{3k+2} x_{3k+1}^{-1}
 x_{3k-1 }^{-1} \ldots  x_{k+3}^{-1} x_{k+1}^{-1} x_0^{-k}$
 of swing $k$. \label{simplerseesaw}}
\end{figure}

The seesaw words used in the proof of Theorem \ref{thm:seesaw} are potentially asymmetric, in that
the two parameters $l$ and $m$ separately control the extent to which
$x_0^{-1}$ and $x_0$ respectively reduce word length.  For simplicity, we
can consider a one-parameter family of seesaw words
 where
we set $l=m=k$ to get seesaw words of swing $k$ in both
directions. These words are of the form $x_0^{k-1} x_1 x_{3k+2}
x_{3k+1}^{-1} x_{3k-1 }^{-1} \ldots  x_{k+3}^{-1} x_{k+1}^{-1}
x_0^{-k}$ and are pictured in Figure \ref{simplerseesaw}.

The existence of these seesaw words eliminates the possibility of
families of geodesics in the Cayley graph $\Gamma(F, \{x_0,x_1\})$
which satisfy the $k$-fellow traveller property.

\begin{proposition}
\label{prop:notkft} Given any constant $k$, there is $w \in \s$ so
that $w, w x_0$ and $w x_0^{-1}$ cannot be represented by geodesic
paths from the identity which satisfy the $k$-fellow traveller
property.
\end{proposition}

\begin{proof} Let $w \in \s$ be a seesaw word of swing $m$.
 Then any geodesic path from the identity to $w$
in the Cayley graph $\Gamma$ must end either in the suffix $x_0^m$
or $x_0^{-m}$.  The first of these possible suffixes comes from a path
which passes through $wx_0^{-1}$ and the second from a path
through $wx_0$.

Let $\gamma$ be a geodesic path from the identity to $w$ passing
through $wx_0$ and $\eta'$ a path from the identity to $w$ passing
though $wx_0^{-1}$. Let $\eta$ be the prefix of this path ending
at $wx_0^{-1}$.  Then the length of $\gamma$ is one more than the
length of $\eta$, and $d_{\Gamma}(w,wx_0^{-1}) = 1$.

We can write $\gamma = \gamma_1 x_0^s$ and $\eta = \eta_1
x_0^{-(s-1)}$, where $s \leq m$. Then we know that $\gamma_1$ and
$\eta_1$ have the same length. To compute the distance
$D_s(\gamma,\eta)$, we consider $d(wx_0^{-s},wx_0^s) = 2s$.  Since
we can find seesaw words of arbitrarily large swing, we can make
this distance arbitrarily large.  Thus $D_s(\gamma,\eta)$ is not
bounded by a constant, and the paths cannot satisfy the $k$-fellow
traveller property for the given constant $k$.
\end{proof}

The paths discussed above show explicitly that $F$ does not admit
a combing by geodesics.  This also follows from the fact that $F$
is not almost convex \cite{ct}.

\begin{theorem}
Thompson's group $F$ is not combable by geodesics.
\end{theorem}

\begin{proof}
Consider any combing of the Cayley graph $\Gamma$ by geodesics.
Let $w \in \s$, and $\gamma$ the geodesic combing path from the
identity to $w$ in $\Gamma$.  Then $\gamma$ passes through $wx_0$
or $wx_0^{-1}$ but not both.  Let $\eta$ be the combing path to
the point $wx_0^{\pm 1}$ not on $\gamma$.  It follows from
Proposition \ref{prop:notkft} that these paths do not satisfy the
$k$ fellow traveller property, and thus $\Gamma$ is not combable
by geodesics.
\end{proof}

%\bibliography{thompc}

\begin{thebibliography}{10}

\bibitem{belkbrown}
James Belk and Kenneth~S. Brown.
\newblock Forest diagrams for elements of {T}hompson's group ${F}$,
  arXiv:math.GR/0305412.
\newblock Preprint.

\bibitem{belkbux}
James Belk and Kai-Uwe Bux.
\newblock Thompson's group ${F}$ is not minimally almost convex,
  arXiv:math.GR/0301141.
\newblock Preprint.

\bibitem{bg:thomp}
Kenneth~S. Brown and Ross Geoghegan.
\newblock An infinite-dimensional torsion-free ${FP_\infty}$ group.
\newblock {\em Inventiones mathematicae}, 77:367--381, 1984.

\bibitem{josegrowth}
Jos{\'e} Burillo.
\newblock Growth of positive words in {T}hompson's group ${F}$.
\newblock Preprint.

\bibitem{cfp}
James~W. Cannon, William~J. Floyd, and Walter~R. Parry.
\newblock Introductory notes on {R}ichard {T}hompson's groups.
\newblock {\em L'Ens. Math.}, 42:215--256, 1996.

\bibitem{ct}
Sean Cleary and Jennifer Taback.
\newblock Thompson's group {$F$} is not almost convex.
\newblock {\em J. Algebra}, 270(1):133--149, 2003.

\bibitem{ctcomb}
Sean Cleary and Jennifer Taback.
\newblock Combinatorial properties of {T}hompson's group {$F$}.
\newblock {\em Trans. Amer. Math. Soc.}, 356(7):2825--2849 (electronic), 2004.

\bibitem{ctdeadlamp}
Sean Cleary and Jennifer Taback.
\newblock {Dead end words in lamplighter groups and other wreath products}.
\newblock {\em Quarterly Journal of Mathematics}, to appear,
  arXiv:math.GR/0309344.

\bibitem{blake}
Blake Fordham.
\newblock {\em Minimal Length Elements of {T}hompson's group ${F}$}.
\newblock PhD thesis, Brigham Young Univ, 1995.

\bibitem{blakegd}
S.~Blake Fordham.
\newblock Minimal length elements of {T}hompson's group {$F$}.
\newblock {\em Geom. Dedicata}, 99:179--220, 2003.

\bibitem{gubaquad}
Victor Guba.
\newblock The {D}ehn function of {T}hompson's group $f$ is quadratic,
  arXiv:math.GR/0211395.
\newblock Preprint.

\bibitem{gubagrowth}
Victor Guba.
\newblock On the properties of the {C}ayley graph of {R}ichard {T}hompson's
  group {F}, arXiv:math.GR/0211396.
\newblock Preprint.

\end{thebibliography}
\bibliographystyle{hplain}

\end{document}